\documentclass[12pt,twoside]{article}

\usepackage[english,russian]{babel}
\usepackage{umj_1}           

\usepackage{graphicx, srcltx}   
\usepackage{cite}               

\usepackage{amsmath}            
\usepackage{amsfonts,amssymb}   
\usepackage{srcltx}             

\newcommand{\Leftclt}{Roman A. Veprintsev}
\newcommand{\Rightclt}{On an inequality of different metrics for algebraic polynomials }
\usepackage{authblk}

\DeclareMathOperator*{\esssup}{ess\,sup}

\begin{document}

\hfill \textsc{In memory of my Grandmother}

\begin{Titul}
{\large \bf ON AN INEQUALITY OF DIFFERENT METRICS\\[2mm] FOR ALGEBRAIC POLYNOMIALS }\\[3ex]
{{\bf Roman~A.~Veprintsev} \\[5ex]}
\end{Titul}

\begin{Anot}
{\bf Abstract.} We establish an inequality of different metrics for algebraic polynomials.


\end{Anot}


%

\section{Introduction and preliminaries}

In this section, we give some notation used in the article.

Consider the generalized Jacobi weight
\begin{equation*}
\omega_{\alpha,\beta,\gamma}(x)=(1-x)^\alpha(1+x)^\beta|x|^\gamma,\qquad x\in[-1,1],
\end{equation*}
where $\alpha,\,\beta,\,\gamma>-1$. Given $1\leq p\leq\infty$, we denote by $L_p(\omega_{\alpha,\beta,\gamma})$ the space of complex-valued Lebesgue measurable functions $f$ on $[-1,1]$ with finite norm
\begin{equation*}
\begin{array}{lr}
\|f\|_{L_p(\omega_{\alpha,\beta,\gamma})}=\Bigl(\int\nolimits_{-1}^1 |f(x)|^p\,\omega_{\alpha,\beta,\gamma}(x)\,dx\Bigr)^{1/p},&\quad 1\leq p<\infty,\\[1.0em]
\|f\|_{L_\infty(\omega_{\alpha,\beta,\gamma})}=\esssup\limits_{x\in[-1,1]} |f(x)|.&
\end{array}
\end{equation*}

Define the uniform norm of a continuous function $f$ on $[-1,1]$ by
\begin{equation*}
\|f\|_{\infty}=\max\limits_{-1\leq x\leq 1} |f(x)|.
\end{equation*}
The maximum and the mimimum of two real numbers $x$ and $y$ are denoted by $\max(x,y)$ and $\min(x,y)$, respectively. For $\alpha,\,\mu\geq0$, $p\in[1,\infty)$, and $n=1,\,2,\,\ldots$, let
\begin{equation*}
l_{\alpha,\mu}=\frac{\alpha}{\alpha+\mu}\quad (l_{0,0}=0),\qquad l_{\alpha,\mu}^{\max}=\frac{\max(\alpha,\mu)}{\alpha+\mu}\quad (l_{0,0}^{\max}=0),
\end{equation*}
\begin{equation}\label{special_constant}
C(\alpha,\mu,p,n)=\Bigl(1-\frac{1}{\pi n}\Bigr)^{-\frac{\min(\alpha,\mu)}{p}}\, 2^{\,1+\frac{1}{p}}\,(\max(\alpha,\mu)+1)^{\frac{1}{p}}\,\pi^{\frac{\max(\alpha,\mu)}{p}}.
\end{equation}
Note that $l_{\alpha,\mu}^{\max}\in\bigl[\frac{1}{2},1\bigr]$ if $\max(\alpha,\mu)>0$.

The aim of the paper is to establish an inequality of different metrics for algebraic polynomials. In order to realize this aim, we prove a generalization of a lemma by N.\,K.~Bari.

\section{Auxiliary results}

In this section, we establish some lemmas that will be used to prove our main results.

\begin{lemen}\label{first_lemma_for_segment}
Let $\alpha\geq\mu\geq0$. Suppose that $\Delta\subset[0,1]$ is any segment of length $l$ with $l\leq l_{\alpha,\mu}$. Then the following inequality holds:
\begin{equation*}
\int\nolimits_0^{l} x^\alpha(1-x)^\mu\,dx\leq \int\nolimits_{\Delta}x^\alpha(1-x)^\mu\,dx.
\end{equation*}
\end{lemen}

\proofen The claim is obviously true when $\mu=0$. We now prove the claim for $\mu>0$.

Note that $l_{\alpha,\mu}\in\bigl[\frac{1}{2},1\bigr)$ when $\mu>0$.
The function $x^\alpha(1-x)^\mu$ is increasing on $\bigl[0,\frac{\alpha}{\alpha+\mu}\bigr]$ and is decreasing on $\bigl[\frac{\alpha}{\alpha+\mu},1\bigr]$, because the derivative of this function is positive on $\bigl(0,\frac{\alpha}{\alpha+\mu}\bigr)$ and is negative on $\bigl(\frac{\alpha}{\alpha+\mu},1\bigr)$.

Let us consider the following cases:
\begin{itemize}
\item[I)] $\Delta\subset[0,l_{\alpha,\mu}]$;

\item[II)] $\Delta\subset[l_{\alpha,\mu},1]$;

\item[III)] $\Delta=[a,b]$ and $l_{\alpha,\mu}\in(a,b)$.
\end{itemize}

Case I). Because the function $x^\alpha(1-x)^\mu$ is increasing on $\bigl[0,\frac{\alpha}{\alpha+\mu}\bigr]$, we have, for any segments $\Delta_1=[a_1,b_1]\subset\bigl[0,\frac{\alpha}{\alpha+\mu}\bigr]$ and $\Delta_2=[a_2,b_2]\subset\bigl[0,\frac{\alpha}{\alpha+\mu}\bigr]$ of equal length with $a_1\leq a_2$ (or, equivalently, $b_1\leq b_2$),
\begin{equation}\label{equation_for_first_case}
\int\nolimits_{\Delta_1} x^\alpha(1-x)^\mu\,dx\leq \int\nolimits_{\Delta_2} x^\alpha(1-x)^\mu\,dx.
\end{equation}
Putting $\Delta_1=[0,l]$ and $\Delta_2=\Delta$ in \eqref{equation_for_first_case}, we obtain the desired inequality.

Case II). Note that in this case $l\leq 1-l_{\alpha,\mu}\leq\frac{1}{2}$. Because the function $x^\alpha(1-x)^\mu$ is decreasing on $\bigl[\frac{\alpha}{\alpha+\mu},1\bigr]$,
\begin{equation}\label{first_equation_for_second_case}
\int\nolimits_{\Delta} x^\alpha(1-x)^\mu\,dx\geq\int\nolimits_{1-l}^1 x^\alpha(1-x)^\mu\,dx.
\end{equation}
Since $x^\alpha(1-x)^\mu\leq (1-x)^\alpha x^\mu$ on $\bigl[0,\frac{1}{2}\bigr]$, we get
\begin{equation}\label{second_equation_for_second_case}
\int\nolimits_{1-l}^1 x^\alpha(1-x)^\mu\,dx=\int\nolimits_0^l (1-x)^\alpha x^\mu\,dx\geq \int\nolimits_{0}^l x^\alpha(1-x)^\mu\,dx.
\end{equation}
Now the desired inequality follows from the inequalities \eqref{first_equation_for_second_case} and \eqref{second_equation_for_second_case}.

Case III). We have
\begin{equation}\label{first_equation_for_third_case}
\int\nolimits_0^l x^\alpha(1-x)^\mu\,dx=\int\nolimits_{0}^{b-l_{\alpha,\mu}} x^\alpha(1-x)^\mu\,dx+\int\nolimits_{b-l_{\alpha,\mu}}^l x^\alpha(1-x)^\mu\,dx
\end{equation}
and
\begin{equation}\label{second_equation_for_third_case}
\int\nolimits_{\Delta} x^\alpha(1-x)^\mu\,dx=\int\nolimits_a^{l_{\alpha,\mu}} x^\alpha(1-x)^\mu\,dx+\int\nolimits_{l_{\alpha,\mu}}^b x^\alpha(1-x)^\mu\,dx.
\end{equation}
Applying \eqref{equation_for_first_case}--\eqref{second_equation_for_second_case} in the appropriate settings, we can obtain
\begin{equation*}
\int\nolimits_{b-l_{\alpha,\mu}}^l x^\alpha(1-x)^\mu\,dx\leq\int\nolimits_a^{l_{\alpha,\mu}} x^\alpha(1-x)^\mu\,dx,
\end{equation*}
\begin{equation*}
\int\nolimits_{0}^{b-l_{\alpha,\mu}} x^\alpha(1-x)^\mu\,dx\leq \int\nolimits_{l_{\alpha,\mu}}^b x^\alpha(1-x)^\mu\,dx.
\end{equation*}
Using \eqref{first_equation_for_third_case}, \eqref{second_equation_for_third_case}, and the above inequalities, we get the desired inequality.
\hfill$\square$

\begin{coren}\label{first_corollary_for_segment}
Let $\mu\geq\alpha\geq0$. Suppose that $\Delta\subset[0,1]$ is any segment of length $l$ with $l\leq l_{\mu,\alpha}$. Then the following inequality holds:
\begin{equation*}
\int\nolimits_{1-l}^{1} x^\alpha(1-x)^\mu\,dx\leq \int\nolimits_{\Delta}x^\alpha(1-x)^\mu\,dx.
\end{equation*}
\end{coren}

\begin{coren}\label{second_corollary_for_segment}
Let $\alpha,\,\mu\geq0$. Suppose that $\Delta\subset[0,1]$ is any segment of length $l$ with $l\leq l_{\alpha,\mu}^{\max}$ and $l<1$.
Then the following inequality holds:
\begin{equation*}
\int\nolimits_{\Delta} x^\alpha(1-x)^\mu\,dx\geq (1-l)^{\min(\alpha,\mu)}\frac{l^{\max(\alpha,\mu)+1}}{\max(\alpha,\mu)+1}.
\end{equation*}
\end{coren}

\proofen If $\alpha\geq\mu$, then, by Lemma~\ref{first_lemma_for_segment}, we get
\begin{equation}\label{first_estimate_for_second_corollary}
\int\nolimits_{\Delta}x^\alpha(1-x)^\mu\,dx\geq \int\nolimits_0^{l} x^\alpha(1-x)^\mu\,dx\geq (1-l)^\mu\int\nolimits_0^l x^\alpha\,dx=(1-l)^\mu \frac{l^{\alpha+1}}{\alpha+1}.
\end{equation}

If $\mu\geq\alpha$, then, by Corollary~\ref{first_corollary_for_segment}, we get
\begin{equation}\label{second_estimate_for_second_corollary}
\int\nolimits_{\Delta}x^\alpha(1-x)^\mu\,dx\geq \int\nolimits_{1-l}^{1} x^\alpha(1-x)^\mu\,dx\geq (1-l)^\alpha\int\nolimits_{1-l}^1 (1-x)^\mu\,dx=(1-l)^\alpha \frac{l^{\mu+1}}{\mu+1}.
\end{equation}

Combining \eqref{first_estimate_for_second_corollary} and \eqref{second_estimate_for_second_corollary}, we obtain the desired estimate.
\hfill$\square$

\begin{lemen}\label{additional_lemma}Let $\alpha,\,\mu\geq0$. Suppose that $\Delta\subset\bigl[0,\frac{\pi}{2}\bigr]$ is any segment of length $l$ with $l\leq \frac{\pi}{2}\,l_{\alpha,\mu}^{\max}$ and $l<\frac{\pi}{2}$. Then
\begin{equation*}
\int\nolimits_{\Delta} |\sin t|^\alpha|\cos t|^\mu\,dt\geq  \Bigl(1-\frac{2l}{\pi}\Bigr)^{\min(\alpha,\mu)} \frac{2^{\,\max(\alpha,\mu)}}{\pi^{\max(\alpha,\mu)}}\,\cdot\,\frac{l^{\max(\alpha,\mu)+1}}{\max(\alpha,\mu)+1}.
\end{equation*}
\end{lemen}

\proofen It is well known that, for $t\in\bigl[0,\frac{\pi}{2}\bigr]$,
\begin{equation*}
\sin t\geq \frac{2t}{\pi},\qquad
\cos t\geq 1-\frac{2t}{\pi}.
\end{equation*}

Let $\Delta=[a,b]$. Note that $b-a=l$, $\bigl(\frac{2b}{\pi}-\frac{2a}{\pi}\bigr)\leq l_{\alpha,\mu}^{\max}$, and $\bigl(\frac{2b}{\pi}-\frac{2a}{\pi}\bigr)<1$. Using the above inequalities and Corollary~\ref{second_corollary_for_segment}, we obtain
\begin{equation*}
\begin{split}
\int_a^b |\sin t|^\alpha|\cos t|^\mu\,dt&\geq \int\nolimits_a^b \Bigl(\frac{2}{\pi}\,t\Bigr)^\alpha \Bigl(1-\frac{2}{\pi}\,t\Bigr)^\mu\,dt=\frac{\pi}{2}\,\int\nolimits_{(2a)/\pi}^{(2b)/\pi} x^\alpha(1-x)^\mu\,dx\geq\\
&\geq \frac{\pi}{2}\,\Bigl(1-\frac{2l}{\pi}\Bigr)^{\min(\alpha,\mu)}\,\cdot\,\frac{\bigl(\frac{2}{\pi}\,l\bigr)^{\max(\alpha,\mu)+1}}{\max(\alpha,\mu)+1}=\\
&=\Bigl(1-\frac{2l}{\pi}\Bigr)^{\min(\alpha,\mu)} \frac{2^{\,\max(\alpha,\mu)}}{\pi^{\max(\alpha,\mu)}}\,\cdot\,\frac{l^{\max(\alpha,\mu)+1}}{\max(\alpha,\mu)+1}.
\end{split}
\end{equation*}
\hfill$\square$

\section{Main results}

The following lemma generalizes Lemma~1 in \cite{bari_article_generalization_of_inequalities_1954}.

\begin{lemen}\label{generalization_of_Bari}Let $\alpha,\,\mu\geq0$, $p\geq1$, $n$ is a positive integer. For any trigonometric polynomial $T_n$ of degree $n$, we have
\begin{equation*}
\max\limits_{-\pi\leq t\leq\pi} |T_n(t)|\leq C(\alpha,\mu,p,n)\,n^{\frac{\max(\alpha,\mu)+1}{p}}\,\Bigl(\int\nolimits_{-\pi}^{\pi} |T_n(t)|^p |\sin t|^{\alpha}|\cos t|^\mu\,dt\Bigr)^{1/p},
\end{equation*}
where the constant $C(\alpha,\mu,p,n)$ is defined in \eqref{special_constant}.
\end{lemen}

\proofen Let
\begin{equation}\label{eq_1_Bari}
\nu=|T_n(t_0)|=\max\limits_{-\pi\leq t\leq\pi} |T_n(t)|,
\end{equation}
\begin{equation*}
\Delta_0=\Bigl[t_0-\frac{1}{2n},t_0+\frac{1}{2n}\Bigr].
\end{equation*}

From Bernstein's inequality it follows that
\begin{equation}\label{eq_2_Bari}
|T'_n(t)|\leq n\nu,\quad t\in[-\pi,\pi].
\end{equation}

It is known that, for any $h\geq0$, there exists a $\theta\in(0,1)$ such that
\begin{equation}\label{eq_3_Bari}
\bigl||T_n(t_0+h)|-|T_n(t_0)|\bigr|\leq|T_n(t_0+h)-T_n(t_0)|=|h||T'_n(t_0+\theta h)|.
\end{equation}
Using \eqref{eq_1_Bari}--\eqref{eq_3_Bari}, we get $\bigl||T_n(t_0+h)|-\nu\bigr|\leq n\nu|h|$. Hence, for $|h|\leq\frac{1}{2n}$,
\begin{equation*}
|T_n(t)|\geq\frac{\nu}{2},\qquad t\in\Delta_0.
\end{equation*}
Thus, we have
\begin{equation}\label{eq_4_Bari}
\begin{split}
\int\nolimits_{-\pi}^{\pi} |T_n(t)|^p |\sin t|^{\alpha}|\cos t|^\mu\,dt&\geq\int\nolimits_{\Delta_0} |T_n(t)|^p |\sin t|^{\alpha}|\cos t|^\mu\,dt\geq\\&\geq\Bigl(\frac{\nu}{2}\Bigr)^p \int\nolimits_{\Delta_0} |\sin t|^{\alpha}|\cos t|^\mu\,dt.
\end{split}
\end{equation}

Since $|\sin t|^{\alpha}|\cos t|^\mu$ is an even function of period $\pi$, we can assume, without loss of generality, that the centre of $\Delta_0$ belongs to $\bigl[0,\frac{\pi}{2}\bigr]$. Then there exists a segment $\Delta$ of length $\frac{1}{2n}$ such that $\Delta\subset\bigl[0,\frac{\pi}{2}\bigr]$ and $\Delta\subset\Delta_0$. Note that
$\frac{1}{2n}<\frac{\pi}{4}\leq \frac{\pi}{2} l_{\alpha,\mu}^{\max}$. Hence, using Lemma~\ref{additional_lemma}, we get
\begin{equation}\label{eq_5_Bari}
\begin{split}
\int\nolimits_{\Delta_0} |\sin t|^{\alpha}|\cos t|^{\mu}\,dt&\geq \int\nolimits_{\Delta} |\sin t|^{\alpha}|\cos t|^{\mu}\,dt\geq\\
&\geq \Bigl(1-\frac{1}{\pi n}\Bigr)^{\min(\alpha,\mu)}\,\cdot\,\frac{1}{2(\max(\alpha,\mu)+1) \pi^{\max(\alpha,\mu)} n^{\max(\alpha,\mu)+1}}.
\end{split}
\end{equation}

From \eqref{eq_4_Bari}, \eqref{eq_5_Bari} it follows that
\begin{equation*}
\begin{split}
\Bigl(\int\nolimits_{-\pi}^{\pi} &|T_n(t)|^p |\sin t|^{\alpha}|\cos t|^\mu\,dt\Bigr)^{1/p}\geq\\
&\geq\frac{\nu}{2}\,\cdot\,\Bigl\{\Bigl(1-\frac{1}{\pi n}\Bigr)^{\min(\alpha,\mu)}\,\cdot\,\frac{1}{2(\max(\alpha,\mu)+1) \pi^{\max(\alpha,\mu)} n^{\max(\alpha,\mu)+1}}\Bigr\}^{1/p}=\\
&=\Bigl\{\Bigl(1-\frac{1}{\pi n}\Bigr)^{-\frac{\min(\alpha,\mu)}{p}} 2^{\,1+\frac{1}{p}}(\max(\alpha,\mu)+1)^{\frac{1}{p}}\pi^{\frac{\max(\alpha,\mu)}{p}} n^{\frac{\max(\alpha,\mu)+1}{p}}\Bigr\}^{-1}\,\cdot\,\nu.
\end{split}
\end{equation*}
\hfill$\square$

Now we list some properties of $C(\alpha,\mu,p,n)$:
\begin{itemize}
\item[$(1)$] $C(\alpha,\mu,p,n)\leq C(\alpha,\mu,p,1)$, $n=1,\,2,\,\ldots$.

\item[$(2)$] $C(\alpha,\mu,p,n)\to 2^{1+\frac{1}{p}}\,(\max(\alpha,\mu)+1)^{\frac{1}{p}}\,\pi^{\frac{\max(\alpha,\mu)}{p}}$, $n\to\infty$.

\item[$(3)$] If $\max(\alpha,\mu)\leq p$, then $C(\alpha,\mu,p,1)\leq\frac{8\pi^2}{\pi-1}$ and $C(\alpha,0,p,1)\leq8\pi$.
\end{itemize}

The following theorem generalizes Lemma (an inequality of different metrics for polynomials) in \cite{kamzolov_article_Foirier-Jacobi_series_2007}.

\begin{teoen}\label{main_theorem}Let $\alpha\geq\beta\geq-\frac{1}{2}$, $\mu\geq0$, $1\leq p<q\leq\infty$, $n$ is a positive integer.
If $P_n$ is an algebraic polynomial of degree $n$, then
\begin{equation*}
\|P_n\|_{L_q(\omega_{\alpha,\beta,\mu})}\leq B(\alpha,\beta,\mu,p,n)^{\left(\frac{1}{p}-\frac{1}{q}\right)} n^{\max(2(\alpha+1),\mu+1)\left(\frac{1}{p}-\frac{1}{q}\right)} \|P_n\|_{L_p(\omega_{\alpha,\beta,\mu})},
\end{equation*}
where
\begin{equation*}
B(\alpha,\beta,\mu,p,n)=2^{\,2p+1+\alpha-\beta} \Bigl(1-\frac{1}{\pi n}\Bigr)^{-\min(2\alpha+1,\mu)} \max(2(\alpha+1),\mu+1) \, \pi^{\max(2\alpha+1,\mu)}.
\end{equation*}
\end{teoen}

\proofen Note that
\begin{equation*}
B(\alpha,\beta,\mu,p,n)=\Bigl\{2^{1+\frac{\alpha-\beta}{p}} C(2\alpha+1,\mu,p,n)\Bigr\}^p.
\end{equation*}

Using Lemma~\ref{generalization_of_Bari}, we get
\begin{equation*}
\begin{split}
\|P_n\|_\infty&=\max\limits_{-\pi\leq t\leq\pi} |P_n(\cos t)|\leq\\&\leq C(2\alpha+1,\mu,p,n) \, n^{\frac{\max(2\alpha+1,\mu)+1}{p}} \Bigl(\int\nolimits_{-\pi}^{\pi} |P_n(\cos t)|^p |\sin t|^{2\alpha+1}|\cos t|^{\mu}\,dt\Bigr)^{1/p}=\\
&=2\, C(2\alpha+1,\mu,p,n) \, n^{\frac{\max\left(2(\alpha+1),\mu+1\right)}{p}}\, \Bigl(\int\nolimits_{0}^{\pi} |P_n(\cos t)|^p (\sin t)^{2\alpha+1}|\cos t|^{\mu}\,dt\Bigr)^{1/p}=\\
&=2^{1+\frac{2\alpha+1}{p}} \, C(2\alpha+1,\mu,p,n) \, n^{\frac{\max\left(2(\alpha+1),\mu+1\right)}{p}} \times\\&\quad\times \Bigl(\int\nolimits_{0}^{\pi} |P_n(\cos t)|^p \Bigl(\sin \frac{t}{2}\Bigr)^{2\alpha+1}\Bigl(\cos \frac{t}{2}\Bigr)^{2\alpha+1}|\cos t|^{\mu}\,dt\Bigr)^{1/p}=\\
&=2 \, C(2\alpha+1,\mu,p,n) \, n^{\frac{\max\left(2(\alpha+1),\mu+1\right)}{p}}\times\\&\quad\times  \Bigl(\int\nolimits_{0}^{\pi} |P_n(\cos t)|^p (1-\cos t)^{\alpha+1/2}(1+\cos t)^{\alpha+1/2}|\cos t|^{\mu}\,dt\Bigr)^{1/p}=\\
&=2^{\,1+\frac{\alpha-\beta}{p}} \, C(2\alpha+1,\mu,p,n) \, n^{\frac{\max\left(2(\alpha+1),\mu+1\right)}{p}}\times\\&\quad\times  \Bigl(\int\nolimits_{0}^{\pi} |P_n(\cos t)|^p (1-\cos t)^{\alpha+1/2}(1+\cos t)^{\beta+1/2}|\cos t|^{\mu}\,dt\Bigr)^{1/p}=\\
&=2^{\,1+\frac{\alpha-\beta}{p}} \, C(2\alpha+1,\mu,p,n) \, n^{\frac{\max\left(2(\alpha+1),\mu+1\right)}{p}}\times\\&\quad\times
\Bigl(\int\nolimits_{-1}^1 |P_n(x)|^p\,(1-x)^{\alpha}(1+x)^{\beta}|x|^{\mu}\,dx\Bigr)^{1/p}=\\
&=2^{\,1+\frac{\alpha-\beta}{p}} \, C(2\alpha+1,\mu,p,n) \, n^{\frac{\max\left(2(\alpha+1),\mu+1\right)}{p}}\,\|P_n\|_{L_p(\omega_{\alpha,\beta,\mu})}.
\end{split}
\end{equation*}
Consequently,
\begin{equation*}
\begin{split}
\|P_n\|_{L_q(\omega_{\alpha,\beta,\mu})}&=\Bigl(\int\nolimits_{-1}^1 |P_n(x)|^q\,(1-x)^{\alpha}(1+x)^{\beta}|x|^{\mu}\,dx\Bigr)^{1/q}\leq\\
&\leq \Bigl(\int\nolimits_{-1}^1 \|P_n\|_{\infty}^{q-p} |P_n(x)|^p\,(1-x)^{\alpha}(1+x)^{\beta}|x|^{\mu}\,dx\Bigr)^{1/q}=\\
&=\|P_n\|_{\infty}^{1-\frac{p}{q}} \|P_n\|_{L_p(\omega_{\alpha,\beta,\mu})}^{\frac{p}{q}}\leq\\
&\leq \Bigl\{\Bigl(2^{\,1+\frac{\alpha-\beta}{p}} \, C(2\alpha+1,\mu,p,n)\Bigr)^p\Bigr\}^{\frac{1}{p}-\frac{1}{q}} n^{\max(2(\alpha+1),\mu+1)\left(\frac{1}{p}-\frac{1}{q}\right)} \|P_n\|_{L_p(\omega_{\alpha,\beta,\mu})}.
\end{split}
\end{equation*}
\hfill$\square$

Now we list some properties of $B(\alpha,\beta,\mu,p,n)$:
\begin{itemize}
\item[$(a)$] $B(\alpha,\beta,\mu,p,n)\leq B(\alpha,\beta,\mu,p,1)$, $n=1,\,2,\,\ldots$.

\item[$(b)$] $B(\alpha,\beta,\mu,p,n)\to2^{\,2p+1+\alpha-\beta}\, \max(2(\alpha+1),\mu+1)\, \pi^{\max(2\alpha+1,\mu)}$, $n\to\infty$.

\item[$(c)$] If $\max(2\alpha+1,\mu)\leq p$, then $B(\alpha,\beta,\mu,p,1)\leq2^{p+\alpha-\beta} \bigl(\frac{8\pi^2}{\pi-1}\bigr)^{p}$.
\end{itemize}

\section{Conclusion}

Our next aim is to prove based on Lemma~2.2 in \cite{fejzullahu_article_2013} that the inequality in Theorem~\ref{main_theorem} is precise in order.

\begin{Biblioen}

\bibitem{bari_article_generalization_of_inequalities_1954}N.\,K.~Bari, A generalization of the inequalities of S.\,N.~Bernstein and A.\,A.~Markov, \textit{Izv. Akad. Nauk SSSR Ser. Mat.} \textbf{18}\,(2)~(1954), 159--176.
    
\bibitem{fejzullahu_article_2013}B.\,Xh.~Fejzullahu, On orthogonal expansions with respect to the generalized Jacobi weight, \textit{Results. Math.} \textbf{63}~(2013), 1177--1193.

\bibitem{kamzolov_article_Foirier-Jacobi_series_2007}A.\,I.~Kamzolov, A norm of partial sums of Fourier\,--\,Jacobi series for functions from $L_p^{(\alpha,\beta)}$, \textit{Moscow University Mathematics Bulletin}, \textbf{62}\,(6)~(2007), 228--236.
    
\end{Biblioen}

\noindent \textsc{Independent researcher, Uzlovaya, Russia}

\noindent \textit{E-mail address}: \textbf{veprintsevroma@gmail.com}

\end{document}